\documentclass{birkjour_t2}
\usepackage{graphicx}
\usepackage{mathrsfs}
\usepackage{amsmath}
\usepackage{color}
\newtheorem{theo}{Theorem}[section]

\newtheorem{lem}[theo]{Lemma}

\newcommand{\R}{\mathbb{R}}
\newcommand{\CC}{\boldsymbol{C}}
\newcommand{\BB}{\boldsymbol{B}}
\newcommand{\HH}{\boldsymbol{H}}
\newcommand{\YY}{\boldsymbol{Y}}
\newcommand{\PPsi}{\boldsymbol{\Psi}}
\newcommand{\TTet}{\boldsymbol{\theta}}
\newcommand{\OM}{\Omega}
\newcommand{\OMV}{\overline{\Omega}}
\newcommand{\FF}{\mathcal{F}}
\newenvironment{pre}{%
	\par\noindent\textbf{Proof. }\ignorespaces
}{%
	\hfill$\square$\par
}
\usepackage{hyperref}

\begin{document}
	
	\title{Existence results for a biofilm free-boundary \\ problem with dominant detachment}
	
	\author{Dieudonn\'e Zirhumanana Balike }
	\address{D\'epartment de Math\'ematique-Physique, Institut Sup\'erieur P\'edagogique de Bukavu, \\30 Avenue Kibombo, 854 Bukavu,Democratic Republic of Congo}
	\email{dieudobalike@gmail.com}
	
	\begin{abstract}
		This work addresses the existence and uniqueness of a Wanner-Gujer free-\\boundary problem that models biofilms under conditions of prevailing detachment. \\This result significantly extends previous findings in both tumor growth modeling and \\the biofilm modeling field.\\
		Besides establishing local existence and uniqueness, we also prove the continuous de-\\pendence of the solution on initial and boundary data. Furthermore, global existence\\ is deduced using a combination of invariance regions and energy estimates. The proof \\ for local existence is obtained by utilizing fixed point arguments combined with\\
		semigroup theory.
	\end{abstract}
	

\maketitle
	\newpage
\section{Introduction}\label{sec1}

We consider the following free bounddary prroblem which is known as Wanner-Gujer model:
\begin{subequations}\label{EqBiomassVersion1}
	\begin{equation}
		\partial_t X_i+u\partial_x X_i =f_i(\boldsymbol{X},\boldsymbol{S}), ~0\le x\le L(t),~t>0,~i=1,\cdots,n; 
	\end{equation}
	\begin{equation}
		X_i(x,0)=\varphi_i(x), ~0\le x\le L_0, ~~i=1,\cdots,n~;
	\end{equation}~
	\begin{equation}
		\partial_x u=g(\boldsymbol{X},\boldsymbol{S}), ~~0<x\le L(t);
	\end{equation}
	\begin{equation}
		u(0,t)=0,~~t >0.
	\end{equation}
\end{subequations}
\begin{subequations}
	\begin{equation}
		\partial_t S_j -D_j \partial_{xx}S_j=h_j(\boldsymbol{X},\boldsymbol{S}),~0<x<L(t),~t>0,~j=1\cdots,m; 
	\end{equation}
	\begin{equation}
		S_j(x,0)=\theta_j(x),~t>0,~j=1\cdots,m; 
	\end{equation}
	\begin{equation}
		\partial_x S_j(0,t)=0; ~~S_j(L(t),t)=\psi_j(t), ~t>0;~j=1,\cdots,m.
	\end{equation}
\end{subequations}
\begin{subequations}\label{ThicknessEq}
	\begin{equation}
		\dot L(t)=u(L(t),t)-\lambda L^2(t),~\lambda>0,~t>0;
	\end{equation}
	\begin{equation}
		L(0)=L_0.
	\end{equation}
\end{subequations}
This model describes the spatio-temporal evolution of biofilm components. 
Each biomass species $X_i$ grows or decays according to $f_i$, and is transported by the velocity field $u$. 
The substrates $S_j$ diffuse and are consumed by the biomass through the terms $h_j$. 
In general, $f_i$ and $h_j$ follow Monod-type kinetics \cite{baeten2018modelling,picioreanu1998mathematical,picioreanu1998new,boyle1999limits}. 
The biofilm thickness $L(t)$ evolves over time as the biofilm grows and detaches, with the detachment rate proportional to $L^2(t)$. The following  notations $\boldsymbol{X}=(X_1\cdots,X_n)$ and $\boldsymbol{S}=(S_1,\cdots,S_m)$ will be used from now on. More details about the derivation and other comments can be found in \cite{wanner1986multispecies,wanner1996mathematical,vafai2010porous,wanner1995modeling,eberl2020challenges}.
The mathematical investigation of its qualitative properties is less expanded compared to the amount of the works applying it to modeling different process in biotechnology and wastewater treatment \cite{balike2024mathematical,bryers2018modeling,van2002mathematical}.\\
When $\lambda=0$ the problem can be treated as the Stefan problem or as tumor growth modeling problems \cite{schatz1969free,sherman1967free,cui2005global,chen2013two,cui2003free}.  In this paper we consider the case where $\lambda\ne 0$ which is general and frequently used in wastewater treatment and biotechnology engineering.  Indeed, while such models are classical in environmental and bioprocess engineering, 
their analytical well-posedness has received comparatively little attention except  in \cite{dacunto2016qualitative,dacunto2021free}  where it was studied by using method of characteristics. But this approach does not allow to address the global existence due to potential blow-up. Hence, this work addresses the question of local and global existence and uniqueness by using a different technique. More precisely, it establishes a complete local well-posedness theory for the coupled hyperbolic--parabolic free-boundary system, formulated after a  nonlinear rescaling to a fixed reference domain. The proof combines maximal $L^p$-regularity for parabolic systems with a fixed-point argument coupling the transport and diffusion equations through the moving boundary ODE.\\
We unify the analysis under minimal regularity assumptions: 
the results hold for all $p>\tfrac{3}{2}$ in one spatial dimension, which 
is sharp for the embeddings needed to control the nonlinear couplings. 
This extends some works from tumor growth modeling that required stronger assumptions for $\lambda=0.$
For example, it was required that $p>5/2$ in \cite{cui2005global}   and $p>5$  in \cite{cui2005analysis} while a classical $C^2$ and Hölder $C^{\alpha+2,\alpha+\frac{1}{2}}$ regularity was respectively obtained  in \cite{cui2001analysis} and \cite{friedman2005analysis}.\\
The approach used yields not only local existence and uniqueness but also continuous dependence on data, positivity preservation of biomass and substrates,  and a continuation criterion that allows extension to global-in-time solutions 
under dissipative (detachment-dominated) conditions. \\
Finally, the analysis proposed  covers both the parabolic and hyperbolic subproblems 
as separate lemmas (Lemmas~\ref{lem:main-final}--\ref{Lemma3.2}), providing reusable well-posedness results 
for related coupled reaction--transport systems beyond biofilm modeling.

These results therefore bridge a gap between applied biofilm models widely used 
in wastewater and bioengineering applications and their mathematical analysis,  offering an adapted existence and uniqueness theory consistent with the 
nonlinear coupling and free-boundary dynamics characteristic of the 
Wanner--Gujer modeling.\\
In the following pages the work is organised as follows: the second section introduces the notations used  and presents a reformuation of  the free boundary problem into an initial value problem. In the third section we present primary lemmans and the main result for local existence. Finally, the last section is about the global existence by using the region invariance annd energy estimates.
\section{Notations and reformualtion of the problem}\label{SectReform}
We introduce the notations which will be used in this work:
Throughout the paper let $p> \tfrac{3}{2}$ be fixed. 
For any open interval $\Omega\subset \mathbb{R}$ and integer $k\ge1$, we denote by 
$W^{k,p}(\Omega)$ the usual Sobolev space, with norm \\
$\|u\|_{W^{k,p}(\Omega)} = \sum_{|\alpha|\le k}\|\partial_x^{\alpha}u\|_{L^p(\Omega)}$.\\
For a positive continuous function $L=L(t)$ defined on $[0,T_1]$ and for some $T_1>0$, we set 
\[
\Omega_L := \{(x,t)\in\mathbb{R}\times(0,T_1): 0<x<L(t)\}
\quad \text{and }\quad \overline{\Omega}_L  \text{ its closure.}
\]
For functions defined on $\Omega_L$, we use the standard parabolic Sobolev space
\[
W^{2,1}_p(\Omega_L)
=\bigl\{u\in L^p(\Omega_L): \partial_t u,\; \partial_x u,\;
\partial_{xx}u\in L^p(\Omega_L)\bigr\},
\]
endowed with the norm
\[
\|u\|_{W^{2,1}_p(\Omega_L)} 
:= \|u\|_{L^p(\Omega_L)} 
+ \|\partial_t u\|_{L^p(\Omega_L)}
+ \|\partial_x u\|_{L^p(\Omega_L)}
+ \|\partial_{xx}u\|_{L^p(\Omega_L)} .
\]
Let $D_p(\Omega)$ denote the trace of $W^{2,1}_p(\Omega\times(0,T_1))$ at time $t=0$:
\[
D_p(\Omega)
= \{\,\theta\in L^p(\Omega):
\exists\,u\in W^{2,1}_p(\Omega\times(0,T_1)),\;
u(\cdot,0)=\theta\,\}.
\]
The trace theorem (see for example \cite{amann2019linear,amann2009anisotropic,agresti2023trace})
implies that for $p>\tfrac{3}{2}$,
\[
W^{2,1}_p(\Omega\times(0,T_1))
\hookrightarrow C([0,T_1];W^{2-2/p,p}(\Omega))
\hookrightarrow C(\overline{\Omega}\times[0,T_1]),
\]
so that time traces and boundary values are well defined and continuous.
The norm in $D_p(\Omega)$ is defined by
\[
\|\theta\|_{D_p(\Omega)}
= \inf\bigl\{
T_1^{-1/p}\|u\|_{W^{2,1}_p(\Omega\times(0,T_1))}:
u\in W^{2,1}_p(\Omega\times(0,T_1)),\ u(\cdot,0)=\theta
\bigr\}.
\]
In particular, if $\theta\in W^{2-2/p,p}(\Omega)$, then $\theta\in D_p(\Omega)$ and
$\|\theta\|_{D_p(\Omega)}\le \|\theta\|_{W^{2-2/p,p}(\Omega)}$. We will consider $\Omega=(0,1)$ and $\overline{\Omega}=[0,1]$ (see section \ref{SectReform}).
Since $p>\tfrac{3}{2}$, we have the continuous embedding 
$W^{2,1}_p((0,1)\times(0,T_1)) \hookrightarrow C([0,1]\times[0,T_1])$,
which is the regularity level used in Lemmas~\ref{lem:main-final}--\ref{Lemma3.2},
and Theorems~\ref{thm:main_existence1}--\ref{thm:global}.\\
Hence all subsequent results are formulated for $p> \tfrac{3}{2}$,
which is the minimal condition ensuring the required embeddings and traces
in one spatial dimension.\\
In order to investigate the existence and uniqueness of the problem \eqref{EqBiomassVersion1}--\eqref{ThicknessEq}
we make the following change of variables :
$$z=\frac{x}{L(t)}, ~\tilde{t}=\int_{0}^{t}\frac{d\tau}{L(\tau)}, ~R(\tilde{t})=L(t), ~ C_j(z,\tilde{t})=S_j(x,t), j=1,\cdots,m;$$ $$Y_i(z,\tilde{t})=X_i(x,t), i=1,\cdots,n;~~v(x,\tilde{t})=L(t)u(x,t).$$
To simplify the notations we omit the tilde on $t,$ hence we get the following intial boundary  value problem
\begin{subequations}\label{EqBiomassVersion2}
	\begin{equation}
		\partial_t Y_i-zv(1,t)\partial_z Y_i =R^2(t)f_i(\boldsymbol{Y},\boldsymbol{C})=:F_i(\boldsymbol{Y},\boldsymbol{C}), 0\le z\le 1,~t>0,~i=1,\cdots,n; 
	\end{equation}
	\begin{equation}
		Y_i(z,0)=\varphi_i(z), ~0\le z\le 1, ~~i=1,\cdots,n~;
	\end{equation}
	\begin{equation}
		v(z,t)=u(z,t)-zu(1,t) \qquad \text{for} \quad 0 \le z\le 1,\quad  t> 0;
	\end{equation}
\end{subequations}
\begin{subequations}
	\begin{equation}\label{ParabolicVersion2}
		\begin{array}{r}
			\partial_t C_j -zv(1,t)\partial_zC_j(z,t)-D_j \partial_{zz}C_j=R^2(t)h_j(\boldsymbol{Y},\boldsymbol{C})=:H_j(\boldsymbol{Y},\boldsymbol{C}),\\
			0<z<1,~t>0,~j=1,\cdots,m; 
		\end{array}
	\end{equation}
	\begin{equation}
		C_j(z,0)=\theta_j(z),~t>0,~j=1,\cdots,m; 
	\end{equation}
	\begin{equation}
		\partial_z C_j(0,t)=0; ~~C_j(1,t)=\psi_j(t), ~t>0;~j=1,\cdots,m.
	\end{equation}
\end{subequations}
\begin{subequations}\label{EqVelocity2}
	\begin{equation}
		v(z,t)=R^2(t)\int_{0}^{z}g(\boldsymbol{Y}(\xi,t),\boldsymbol{C}(\xi,t))d\xi
	\end{equation}
	\begin{equation}
		v(1,t)=R^2(t)\int_{0}^{1}g(\boldsymbol{Y}(\xi,t),\boldsymbol{C}(\xi,t))d\xi
	\end{equation}
\end{subequations}
\begin{subequations}\label{EqThicknessVersion2}
	\begin{equation}
		\dot R(t)=R^2(t)v(1,t)-\lambda R^4(t),~t>0;
	\end{equation}
	\begin{equation}
		R(0)=R_0.
	\end{equation}
\end{subequations}
Define the following vector-valued functions :  $\boldsymbol{C}=(C_1,\cdots,C_m)^T,\boldsymbol{Y}=(Y_1,\cdots,Y_n)^T,$ $\boldsymbol{F}=(F_1,\cdots,F_n)^T,$\\$ \boldsymbol{H}=(H_1,\cdots,H_m)^T, \boldsymbol{B}(t)=B(t)\boldsymbol{I}_m$ where $B(t)=zv(1,t)$ and $\boldsymbol{I}_m$ the identity matrix of order $m.$
This system is equivalent to  equations \eqref{EqBiomassVersion1}--\eqref{ThicknessEq} and will be the focus of the rest of the work. We will first give primary lemmas which will be used in the proof of the main results.
\section{Primary lemmas and local existence}
With the notations introduced in the previous section we   recast the parabolic system \eqref{ParabolicVersion2} into a vector form as follows
\begin{equation}\label{eq:main-problem-final}
	\left\{ \begin{array}{r}
		\partial_t \boldsymbol{C}(z,t)=\partial_{zz} \boldsymbol{C}(z,t)+\boldsymbol{B}(t)\partial_z \boldsymbol{C}(z,t)+\boldsymbol{H}(\boldsymbol{Y}(z,t),\boldsymbol{C}(z,t)),\\~(x,t)\in \Omega \times [0,T_1];\\
		\\
		\boldsymbol{C}(z,0)=\boldsymbol{\theta}(z), ~z\in \overline{\Omega};\\
		\\
		\partial_z\boldsymbol{C}(0,t)=0;~~ \boldsymbol{C}(1,t)~=\boldsymbol{\Psi}(t),~t>0.
	\end{array}
	\right.
\end{equation}
where $\Omega=(0,1)$ and $\overline{\Omega}=[0,1].$\\
We make the following assumptions
\begin{itemize}
	\item[(i)]$\BB \in C^\alpha([0,T_1]; \R^{m \times m})$ for some $\alpha \in (0,1];$
	\item[(ii)]  $\HH: \R^m \times \R^m \to \R^m$ is globally Lipschitz with constant $L_H;$
	\item[(iii)]  $\YY \in L^\infty(\Omega\times (0,T_1))^m;$
	\item[(iv)]  $\TTet \in W^{2-2/p,p}(\Omega)^m$ for some $p \in (\frac{3}{2}, \infty);$
	\item[(v)]  $\PPsi \in W^{1-\frac{1}{2p}, p}(0,T_1)^m$ with $\TTet(1) = \PPsi(0);$
	\item[(vi)] (Second-order compatibility)
	$\partial_{zz} \TTet(1) + \BB(0)\partial_z \theta(1) + \HH(\YY(1,0),\TTet(1)) = \frac{d}{dt}\PPsi(0).$
\end{itemize}
\begin{lem}\label{lem:main-final}
	Under the assumptions (i)--(vi) there exists $T^* \in (0, T_1]$ and a unique solution $\CC \in W^{2,1}_p(\Omega \times (0,T^*))^m$ to \eqref{eq:main-problem-final}. Moreover the following estimate holds
	\begin{equation*}
		\|\CC\|_{W^{2,1}_p(\Omega \times (0,T^*))^m} \le M\left( \|\TTet\|_{W^{2-2/p,p}(\Omega)^m} + \|\HH(\YY,0)\|_{L^p(\Omega \times (0,T^*))^m} \right.
	\end{equation*}
	\begin{equation}\label{eq:main-estimate-final}
		\left.+ \|\PPsi\|_{W^{1-\frac{1}{2p},p}(0,T^*)^m} + 1 \right)
	\end{equation}
	where $M > 0$ depends on $p, T^*, \|\BB\|_{C^\alpha([0,T^*])}$, and $L_H$.
\end{lem}

\begin{pre}
	The proof is divided into several steps, the first of which is the parabolic lifiting of boundary data.
	In fact, since $\PPsi \in W^{1-\frac{1}{2p},p}(0,T_1)^m$, by the parabolic trace theorem (see for example  \cite{ladyzhenskaia1968linear,meyries2014traces}), there exists $\tilde{\PPsi} \in W^{2,1}_p(\Omega \times (0,T_1))^m$ such that:
	\begin{itemize}
		\item $\tilde{\PPsi}(1,t) = \PPsi(t)$ for $t \in (0,T_1);$
		\item $\tilde{\PPsi}(z,0) = \TTet(z)$ near $z=1$ (can be arranged via localization);
		\item $\|\tilde{\PPsi}\|_{W^{2,1}_p(\OM \times (0,T_1))^m} \le c\left( \|\PPsi\|_{W^{1-\frac{1}{2p},p}(0,T_1)^m} + \|\TTet\|_{W^{2-2/p,p}(\OM)^m} \right).$
	\end{itemize}
	These conditions will later allow to apply maximal regularity in  this proof. \\
	The second step is the homogenization and abstract formulation by which we will convert the problem into an abstract evolution equation which later will be used in the next steps.\\
	Let us define $\tilde{\CC}(z,t) = \CC(z,t) - \tilde{\PPsi}(z,t)$. Then $\tilde{\CC}$ satisfies:
	\begin{equation}\label{eq:homogeneous-final}
		\begin{cases} 
			\partial_t \tilde{\CC} = \partial_{zz} \tilde{\CC} + \BB(t)\partial_z \tilde{\CC} + \tilde{\HH}(t), & (z,t) \in \OM \times (0,T_1] \\
			\partial_z \tilde{\CC}(0,t) = 0, \quad \tilde{\CC}(1,t) = 0, & t \in (0,T_1] \\
			\tilde{\CC}(z,0) = \tilde{\CC}_0(z) := \CC_0(z) - \tilde{\PPsi}(z,0), & z \in \OM
		\end{cases}
	\end{equation}
	where
	\[
	\tilde{\HH}(z,t) = \HH(\YY(z,t), \tilde{\CC}(z,t) + \tilde{\PPsi}(z,t)) + \BB(t)\partial_z \tilde{\PPsi} + \partial_{zz} \tilde{\PPsi} - \partial_t \tilde{\PPsi}.
	\]
	On $\mathbb{X} = L^p(\Omega)^m$, we define the following differential operator
	\[
	\mathcal{A}(t)\phi = (\partial_{zz}+ \BB(t)\partial_z)\phi \]
	on the domain \[	D(\mathcal{A}(t)) = D := \{\phi \in W^{2,p}(\Omega)^m : \partial_z\phi(0) = 0, \phi(1) = 0\}.\]
	The problem becomes:
	\begin{equation}\label{eq:abstract-final}
		\begin{cases} 
			\frac{d\tilde{\CC}}{dt}(t) = \mathcal{A}(t)\tilde{\CC}(t) + \FF(t, \tilde{\CC}(t)), & t \in (0,T^*] \\
			\tilde{\CC}(0) = \TTet
		\end{cases}
	\end{equation}
	with $\FF(t,\phi)(z) = \HH(\YY(z,t), \phi(z) + \tilde{\PPsi}(z,t)) + (\BB(t)\partial_z+ \partial_{zz} - \partial_t )\tilde{\PPsi}(z,t)$.
	Let $\mathcal{A}_0\phi = \partial_{zz}\phi$ with domain $D$. On one hand, by \cite[Theorem 2.7, Chap.7]{pazy2012semigroups} and \cite[Chap. 8-9]{lorenzi2021semigroups}, $\mathcal{A}_0$ is strongly elliptic and  generates an analytic semigroup on $\mathbb{X}$.\\
	On the otehr hand, the perturbation $P(t)\phi = \BB(t)(z\partial_z\phi)$ is bounded using the one-dimensional  Gagliardo-Nirenberg inequality \cite{dolbeault2014one,nirenberg1959elliptic}:
	\[
	\|\partial_z\phi\|_{L^p} \le C\|\phi\|_{W^{2,p}}^{1/2}\|\phi\|_{L^p}^{1/2} \le \varepsilon\|\phi\|_{W^{2,p}} + C_\varepsilon\|\phi\|_{L^p}
	\]
	we obtain:
	\[
	\|P(t)\phi\|_{L^p} \le \|\BB\|_{L^\infty}\left(\varepsilon\|\mathcal{A}_0\phi\|_{L^p} + C_\varepsilon\|\phi\|_{L^p}\right)
	\]
	Thus $P(t)$ is $\mathcal{A}_0$-bounded with relative bound 0 and  $\mathcal{A}(t)$ generates an analytic semigroup.
	Since $t \mapsto \BB(t)$ is Hölder continuous and $D$ is constant, it follows that the Acquistapace-Terreni conditions \cite{acquistapace1987unified} are satisfied, yielding an evolution system $\{U(t,s)\}_{0 \le s \le t \le T_1}$.
	To prove the existence of a unique  fixed-point we  investigate the maximal regularity. Let us define the Banach space:
	\[
	\mathbb{E}_{T^*} = W^{1,p}(0,T^*; X) \cap L^p(0,T^*; D)
	\]
	which is isomorphic to $W^{2,1}_p(I \times (0,T^*))^m$ by standard parabolic theory. We also need to   the following closed ball  
	\[
	\mathcal{B}_{T^*} = \left\{ \phi \in \mathbb{E}_{T^*} : \phi(0) = \tilde{\CC}_0, \ \|\phi\|_{\mathbb{E}_{T^*}} \le R^* \right\}
	\]
	for some  $R^* > 0.$
	With these tools at hand we can rewrite the abstract evolution  equation \eqref{eq:abstract-final} into an integal form as follows
	\begin{equation}\label{IntegralEq1}
		\Phi(\phi)(t) = U(t,0)\TTet + \int_0^t U(t,s)\FF(s,\phi(s))\, ds
	\end{equation}
	The last part of the proof is dedicated to show that the right hand side of the integral equation \eqref{IntegralEq1} has a unique fixed point.  The following standard arguments are used: first it is easy to show that \cite{pruss2016moving}:
	
	\[
	\|\Phi(\phi)\|_{\mathbb{E}_{T^*}} \le M_1\left( \|\TTet\|_{W^{2-2/p,p}} + \|\FF(\cdot,\phi)\|_{L^p(0,T^*;\mathbb{X})} \right).
	\]
	
	In particular the  the nonlinear term satisify 
	\begin{align*}
		\|\FF(\cdot,\phi)\|_{L^p(0,T^*;\mathbb{X})} &\le L_H\left( \|\phi\|_{L^p(0,T^*;\mathbb{X})} + \|\tilde{\PPsi}\|_{L^p(0,T^*;\mathbb{X})} \right) \\
		&\quad + \|\HH(\YY,0)\|_{L^p(0,T^*;\mathbb{X})} + C\|\tilde{\PPsi}\|_{W^{2,1}_p}
	\end{align*}
	
	By the embedding $\mathbb{E}_{T^*} \hookrightarrow L^\infty(0,T^*; W^{2-2/p,p}(\Omega))$ \cite{lunardi2012analytic} we have 
	\[
	\|\phi\|_{L^p(0,T^*;\mathbb{X})} \le (T^*)^{1/p} \|\phi\|_{L^\infty(0,T^*;\mathbb{X})} \le C(T^*)^{1/p} \|\phi\|_{\mathbb{E}_{T^*}}
	\]
	Using the assumption (ii) we get the following Lipschitz estimate
	\[
	\|\FF(\cdot,\phi_1) - \FF(\cdot,\phi_2)\|_{L^p(0,T^*;\mathbb{X})} \le L_H C(T^*)^{1/p} \|\phi_1 - \phi_2\|_{\mathbb{E}_{T^*}}.
	\]
	
	We therefore choose $R^*$ large enough and $T^*$ small enough so that:
	\[
	M_1\left( \|\tilde{\CC}_0\|_{W^{2-2/p,p}} + L_H(R^* + \|\tilde{\PPsi}\|) + \|\HH(\YY,0)\| + C\|\tilde{\PPsi}\| \right) \le R^*
	\]
	and
	\[
	M_1 L_H C (T^*)^{1/p} < 1
	\]
	Hence  $\Phi$ is a contraction on $\mathcal{B}_{T^*}$.
	By standard bootstrapping and uniqueness of mild/strong solutions, this solution coincides with the unique classical solution in $W_p^{2,1}$\\
	The fixed point $\tilde{\CC} \in \mathbb{E}_{T^*}$ gives the solution to the homogeneous problem. Returning to $\CC = \tilde{\CC} + \tilde{\PPsi}$, we obtain the solution to the original problem.\\
	The estimate follows from the fixed point property and the linear estimates. The second-order compatibility condition ensures that the solution maintains the $W^{2,1}_p$ regularity up to $t=0$.
\end{pre}
\begin{lem}\label{Lemma3.2}[Well-posedness and estimates for a transport–reaction system]
	Let $v \in C^1(\OMV\times[0,T_1])$ satisfy $\|v\|_{L^\infty} + \|\partial_z v\|_{L^\infty} < \infty$ and $v(0,t)=0$ for all $t\in[0,T_1]$.
	Let $F_i : \mathbb{R}^n\times\mathbb{R}^n\times\OMV\times[0,T_1]\to\mathbb{R}$, $i=1,\dots,n$, be Lipschitz in the first two arguments uniformly on $\OMV\times[0,T_1]$, i.e.
	\[
	|F_i(\mathbf{Y},\mathbf{C}) - F_i(\mathbf{\tilde Y},\mathbf{\tilde C})|
	\le L_F\big(|\mathbf{Y}-\mathbf{\tilde Y}| + |\mathbf{C}-\mathbf{\tilde C}|\big),
	\]
	and bounded:
	\[
	\|F_i\|_\infty := \sup_{(\mathbf{Y},\mathbf{C},z,t)} |F_i(\mathbf{Y},\mathbf{C})| < \infty.
	\]
	Assume $\mathbf{C}\in C(\OMV\times[0,T_1];\mathbb{R}^n)$ is given, and $\varphi_i\in C(\OMV)$ for $i=1,\dots,n$.
	
	Then the system
	\begin{equation}\label{eq:transport_system}
		\left\{
		\begin{array}{ll}
			\partial_t Y_i - v(1,t)\,\partial_z Y_i = F_i(\mathbf{Y}(z,t),\mathbf{C}(z,t),z,t), & 0<z<1,\; t>0,\\[0.3em]
			Y_i(z,0)=\varphi_i(z), & 0\le z\le 1,
		\end{array}
		\right.
	\end{equation}
	admits a unique mild (weak) solution $\mathbf{Y}\in C(\OMV\times[0,T_1];\mathbb{R}^n)$, and
	\begin{equation}\label{eq:weak_estimate}
		\|\mathbf{Y}\|_\infty \le e^{L T_1}\big(\|\boldsymbol{\varphi}\|_\infty + T_1\|\mathbf{F}\|_\infty\big).
	\end{equation}
	If in addition each $\varphi_i\in C^1(\OMV )$ and $F_i$ are $C^1$ in $z$ with bounded derivatives, then the weak solution is classical, i.e
	\[
	\mathbf{Y}\in C^1(\OMV\times[0,T_1];\mathbb{R}^n),
	\]
	and the the following estimate holds:
	\begin{equation}\label{eq:classical_estimate}
		\|\mathbf{Y}\|_\infty + \|\partial_z\mathbf{Y}\|_\infty
		\le e^{L T_1}\big(\|\boldsymbol{\varphi}\|_{C^1} + T_1\|\mathbf{F}\|_\infty\big)
		+ T_1 e^{L T_1}\|\partial_z\mathbf{F}\|_\infty.
	\end{equation}
	Moreover, if $\varphi_i(z)\ge 0$ for all $z$ and $F_i(\mathbf{Y},\mathbf{C})\ge 0$ whenever $\mathbf{Y}\ge 0$, then
	\[
	Y_i(z,t)\ge 0,\qquad 0\le z\le 1,\; t\in[0,T_1],\; i=1,\dots,n.
	\]
\end{lem}

\begin{pre}
	Since $v(1,t)$ depends only on $t$, we define the characteristic curve
	\[
	\frac{dZ}{ds} = -v(1,s), \qquad Z(0;z,t)=z.
	\]
	Because $v(1,\cdot)$ is continuous and bounded, there exists a unique absolutely continuous solution:
	\[
	Z(s;z,t) = z - \int_s^t v(1,\tau)\,d\tau.
	\]
	For fixed $(z,t)$ and each $i$, let $\zeta_i(s) = Y_i(Z(s;z,t),s)$. Then
	\[
	\frac{d\zeta_i}{ds}
	= \partial_t Y_i(Z(s),s) + \partial_z Y_i(Z(s),s)\frac{dZ}{ds}
	= \partial_t Y_i - v(1,s)\partial_z Y_i = F_i(\mathbf{Y}(Z(s),s),\mathbf{C}(Z(s),s)).
	\]
	Integrating from $s=0$ to $s=t$ gives the integral formulation
	\begin{equation}\label{eq:integral_form}
		Y_i(z,t)
		= \varphi_i(Z(0;z,t))
		+ \int_0^t F_i(\mathbf{Y}(Z(s;z,t),s),\mathbf{C}(Z(s;z,t),s))\,ds.
	\end{equation}
	Define $\mathcal{T}:C(\OMV\times[0,T_1];\mathbb{R}^n)\to C(\OMV\times[0,T_1];\mathbb{R}^n)$ by
	\[
	(\mathcal{T}\mathbf{Y})_i(z,t)
	:= \varphi_i(Z(0;z,t))
	+ \int_0^t F_i(\mathbf{Y}(Z(s;z,t),s),\mathbf{C}(Z(s;z,t),s))\,ds.
	\]
	For $\mathbf{Y},\mathbf{\tilde Y}\in C$,
	\[
	|(\mathcal{T}\mathbf{Y})_i - (\mathcal{T}\mathbf{\tilde Y})_i|
	\le \int_0^t L |\mathbf{Y}(Z(s),s)-\mathbf{\tilde Y}(Z(s),s)|\,ds
	\le L t \|\mathbf{Y}-\mathbf{\tilde Y}\|_\infty.
	\]
	Hence $\|\mathcal{T}\mathbf{Y}-\mathcal{T}\mathbf{\tilde Y}\|_\infty \le L T_1 \|\mathbf{Y}-\mathbf{\tilde Y}\|_\infty$.
	For $L T_1 < 1$, $\mathcal{T}$ is a contraction, so by the Banach fixed-point theorem,
	there exists a unique solution $\mathbf{Y}$ to~\eqref{eq:integral_form} on $[0,T_1]$.
	This $\mathbf{Y}$ is continuous and is the unique mild (weak) solution of~\eqref{eq:transport_system}.
	The $L^\infty$ estimate are obtained in the following way. 
	From~\eqref{eq:integral_form}, one easily get
	\[
	|Y_i(z,t)|
	\le \|\varphi_i\|_\infty + \int_0^t |F_i(\mathbf{Y}(Z(s),s),\mathbf{C}(Z(s),s))|\,ds.
	\]
	Since $|F_i(\mathbf{Y},\mathbf{C})|\le L|\mathbf{Y}| + \|F_i\|_\infty$, we have
	\[
	|Y_i(z,t)|
	\le \|\varphi_i\|_\infty + L\int_0^t |\mathbf{Y}(Z(s),s)|\,ds + t\|F_i\|_\infty.
	\]
	Let $\psi(t) = \max_i \sup_{z\in[0,1]} |Y_i(z,t)|$. Then
	\[
	\psi(t)\le \|\boldsymbol{\varphi}\|_\infty + t\|\mathbf{F}\|_\infty + L\int_0^t \psi(s)\,ds.
	\]
	By Gronwall’s inequality,
	\[
	\psi(t)\le e^{L t}\big(\|\boldsymbol{\varphi}\|_\infty + t\|\mathbf{F}\|_\infty\big),
	\]
	and setting $t=T_1$ gives~\eqref{eq:weak_estimate}.
	To ensure the existence of a classical solution we assume 
	now that $\varphi_i\in C^1(\OMV)$ and $F_i$ are $C^1$ in $z$ with bounded derivatives.
	Differentiating~\eqref{eq:integral_form} with respect to $z$ gives
	\[
	\partial_z Y_i(z,t)
	= \varphi_i'\big(Z(0;z,t)\big)
	+ \int_0^t \partial_z F_i(\mathbf{Y}(Z(s),s),\mathbf{C}(Z(s),s))\,ds.
	\]
	Hence $\partial_z Y_i$ exists and is continuous, so $\mathbf{Y}\in C^1$. Moreover,
	\[
	|\partial_z Y_i(z,t)| \le \|\varphi_i'\|_\infty + t\|\partial_z F_i\|_\infty.
	\]
	Combining with the previous $L^\infty$ estimate and using again Gronwall’s lemma yields
	\[
	\|\mathbf{Y}\|_\infty + \|\partial_z\mathbf{Y}\|_\infty
	\le e^{L T_1}\big(\|\boldsymbol{\varphi}\|_{C^1} + T_1\|\mathbf{F}\|_\infty\big)
	+ T_1 e^{L T_1}\|\partial_z\mathbf{F}\|_\infty,
	\]
	which is~\eqref{eq:classical_estimate}.
	If $\varphi_i(z)\ge 0$ and $F_i(\mathbf{Y},\mathbf{C})\ge 0$ for $\mathbf{Y}\ge 0$,
	then from~\eqref{eq:integral_form}
	\[
	Y_i(z,t)
	= \varphi_i(Z(0;z,t))
	+ \int_0^t F_i(\mathbf{Y}(Z(s;z,t),s),\mathbf{C}(Z(s;z,t),s))\,ds \ge 0,
	\]
	so the solution remains nonnegative.
	This completes the proof.
\end{pre}
For the local existence we will need the following assumptions in addition to those admititted in Lemmas \ref{lem:main-final} and \ref{Lemma3.2}.\\
(H1) For the initial and boundary and boundary data we assume the follwing
\begin{itemize}
	\item $\boldsymbol{\varphi}=(\varphi_1,\dots,\varphi_n)\in C([0,1];\mathbb R^n)$;
	\item $\TTet=(\theta_1,\dots,\theta_m)\in W^{2-2/p,p}(\OM;\mathbb R^m)$;
	$\PPsi=(\psi_1,\dots,\psi_m)\in W^{1,p}(0,T_1;\mathbb R^m)$ with the compatibility $\theta(1)=\Psi(0)$.
	\item $R_0>0$ is the initial domain size.
\end{itemize}

(H2) For coefficients and the reaction terms we consider the following  regularity assumptions
\begin{itemize}
	\item $v\in C^1(\OMV\times[0,T_1])$ and $\|v\|_\infty+\|\partial_z v\|_\infty<\infty$. 
	\item The reaction maps $\boldsymbol{F}=(F_1,\dots,F_n)$ and $\boldsymbol{H}=(H_1,\dots,H_m)$ are sufficiently regular: 
	$\boldsymbol{F}:\mathbb R^n\times\mathbb R^m\to\mathbb R^n$ and $\boldsymbol{H}:\mathbb R^n\times\mathbb R^m\to\mathbb R^m$ are globally Lipschitz (in their vector arguments) and bounded on bounded sets. Moreover for the classical-regularity part assume $\boldsymbol{F}$ and $\boldsymbol{H}$ are $C^1$ in their arguments with bounded derivatives on the relevant range.
\end{itemize}
Next we give the 
existence and uniqueness theorem  for the full coupled system studied in the project:
The proof uses the two primary lemmas above  and a contraction
argument on a suitable product space.
\begin{theo}\label{thm:main_existence1}
	For $T_1>0$ and under the assumptions made above 
	there exists $T^\ast\in(0,T_1]$ and a unique solution triple
	\[
	(\YY,\CC, v,R)\quad\text{on }\OMV\times[0,T^\ast]
	\]
	such that
	\begin{itemize}
		\item $\CC\in W^{2,1}_p\big(\OM\times(0,T^\ast)\big)^m$ and satisfies the
		parabolic subsystem in the
		classical sense given by Lemma~\ref{lem:main-final};
		\item $\YY\in C^1(\OMV\times[0,T^\ast])^n$ and satisfies equation \eqref{EqBiomassVersion2};
		\item $v\in C^1(\OM\times [0,T^*])$  and $R\in C^1([0,T^\ast])$  satsify respectively \eqref{EqVelocity2} and \eqref{EqThicknessVersion2}  where $v(1,t)$ is the trace of the velocity at $z=1$ computed from $\YY$.
	\end{itemize}
	Moreover there exist constants $K_1,K_2>0$ depending only on the data norms
	\(\|\varphi\|_{C^1},\)\\
	\(\|\theta\|_{W^{2-2/p,p}},\|\psi\|_{W^{1,p}},R_0\), the
	Lipschitz constants, and $p$, such that for $0\le t\le T^\ast$:
	\begin{equation}\label{eq:main_estimate}
		\begin{aligned}
			&\|\CC\|_{W^{2,1}_p((\OM\times(0,t))} 
			+ \|\YY\|_{C^1(\OMV\times[0,T^*])}
			+ \|R\|_{C^1([0,T^*])}
			\\
			&\qquad \le
			K_1 \exp\!\big(K_2 T^*\big)\Big(1+\|\varphi\|_{C^1} + \|\theta\|_{W^{2-2/p,p}} + \|\psi\|_{W^{1,p}} + R_0\Big).
		\end{aligned}
	\end{equation}
	Finally this unique solution  depends continuously on
	the data, and can be extended as long as the norms controlling the estimates
	remain finite.
\end{theo}
\begin{pre}
	The proof relies again on a fixed point argument on a small time interval. We split the construction into maps solved by the two given lemmas, $v$ and for $R$.
	By fixing  \(T>0\) to be chosen and we  define the Banach space
	\[
	\mathcal X_T:=X_Y\times X_R,\qquad
	X_Y:=\{\YY \in C^1(\OMV\times[0,T])^n:\|\YY\|_{C^1}\le M^*\},
	\]
	\[
	X_R:=\{R\in C^1([0,T]):\|R\|_{C^1}\le M^*,\;R(0)=R_0\},
	\]
	with norm \(\|(\YY,R)\|_{\mathcal X_T}:=\|Y\|_{C^1}+\|R\|_{C^1}\). The constant
	\(M^*>0\) can be chosen large enough (depending on the data) so the map we
	construct maps the closed ball \(\mathcal B_M\subset\mathcal X_T\) into itself.
	By choosing arbitrary \((\YY,R)\in\mathcal B_M\) we recall the parabolic right-hand side and transport coefficient
	\[
	H_j(z,t):=R^2(t) h_j\big(\YY(z,t),\CC(z,t)\big),
	\]
	\[\tilde a(t):=v(1,t)=R^2(t)\int_0^1 g(\YY(\xi,t),\CC(\xi,t))\,d\xi.\]
	Note that at this stage $\CC$ is unknown inside the definition of $H_j$; we
	apply the parabolic lemma \ref{lem:main-final} in its mild nonlinear form: the lemma
	assumes the source depends Lipschitz continuously on the unknown \(\CC\) and gives a unique solution \(\CC\in W^{2,1}_p\) together with
	the estimate (for some constant $C_P$ depending on $p$, data and the norms of
	$Y,R$)
	\begin{equation}\label{est:C_from_YR}
		\|\CC\|_{W^{2,1}_p(\OM\times(0,T))}
		\le C_P\Big(\|\TTet\|_{W^{2-2/p,p}} + \|\boldsymbol{H}(\cdot,\cdot,0)\|_{L^p} + \|\PPsi\|_{W^{1,p}}\Big),
	\end{equation}
	where the dependence on \(\YY,R\) is via \(\|\YY\|_{C^1}\) and \(\|R\|_{C^1}\).
	By applying again the parabolic lemma \ref{lem:main-final} in its mild nonlinear form	and by Sobolev embedding (since \(p>3/2\)) we also obtain
	\begin{equation}\label{est:C_uniform}
		\|\CC\|_{C(\OMV\times[0,T])}\le C_{emb}\,\|\CC\|_{W^{2,1}_p((0,1)\times(0,T))}
	\end{equation}
	for some $C_{emb}>0.$
	Thus the parabolic lemma \ref{lem:main-final} furnishes a mapping
	\[
	\mathcal P:\mathcal B_M\to W^{2,1}_p((0,1)\times(0,T))^m,\quad (\YY,R)\mapsto \CC.
	\]
	Using the $\CC$ just obtained and the given \(\YY\), we define
	\begin{subequations}\label{def:v}
		\begin{equation}
			v(z,t):=R^2(t)\int_0^z g\big(\YY(\xi,t),\CC(\xi,t)\big)\,d\xi,
		\end{equation}
		\begin{equation}
			v(1,t)=R^2(t)\int_0^1 g\big(\YY(\xi,t),\CC(\xi,t)\big)\,d\xi.
		\end{equation}
	\end{subequations}
	Because \(g\) is Lipschitz and bounded and \(\YY,\CC,R\) are continuous, the
	integral defines \(v\in C^1(\OMV\times[0,T])\). Moreover, by the uniform
	bounds on \(\YY,\CC,R\) (obtained from the choice of \(\mathcal B_M\) and
	\eqref{est:C_uniform}), the map
	\((\YY,R)\mapsto v\) is Lipschitz into \(C^1\) with a Lipschitz constant that
	is proportional to \(T\) (because of the time integral structure when
	comparing two such \(v\)'s).\\
	With the computed \(v(1,t)\) and the right-hand side
	\[
	F_i(z,t):=R^2(t) f_i\big(\YY(z,t),\CC(z,t)\big),
	\]
	apply the hyperbolic Lemma \ref{Lemma3.2} to get the unique classical solution
	\(\widetilde Y\in C^1(\OMV\times[0,T])^n\) of
	\[
	\partial_t \widetilde Y_i - z v(1,t) \partial_z \widetilde Y_i = F_i(\cdot),
	\qquad \widetilde Y_i(\cdot,0)=\varphi_i.
	\]
	From Lemma \ref{Lemma3.2}  one can prove that there exists  some constant $C_H$ depending on the
	Lipschitz constant of \(f\), the data and the bound \(M\) such that:
	\begin{equation}\label{est:Y_from_C}
		\|\widetilde Y\|_{C^1(\OMV\times[0,T])}
		\le C_H \exp(C_H T)\Big(\|\varphi\|_{C^1} + T\|F\|_{C(\OMV\times[0,T])}\Big).
	\end{equation}
	With these estimates above at hand, we now update $R$ by the ordinary differential equation. Thus, we 
	define \(\widetilde v(1,t)\) from \eqref{def:v} and solve the ODE
	\[
	\dot{\widetilde R}(t)=\widetilde R^2(t)\,\widetilde v(1,t)-\lambda \widetilde R^4(t),
	\qquad \widetilde R(0)=R_0.
	\]
	The right-hand side is locally Lipschitz in \(\widetilde R\) and continuous in
	\(t\), so the ODE has a unique \(C^1\)-solution on \([0,T]\). Moreover,
	using boundedness of \(\widetilde v(1,t)\) and standard ODE estimates we get
	\begin{equation}\label{est:R_from_v}
		\|\widetilde R\|_{C^1([0,T])}\le C_R\big(1+\|v(1,\cdot)\|_{C([0,T])}\big),
	\end{equation}
	with \(C_R\) depending on \(R_0,\lambda\) and \(T\).\\
	We are now in position to define a solution operator that takes into account  all the estimates established so far to deduce the fixed point result. \\ Let us introduce the operator 
	\[
	\mathcal G:\mathcal B_M\to\mathcal X_T,\qquad
	\mathcal G(\YY,R):=(\widetilde \YY,\widetilde R),
	\]
	where \(\widetilde \YY\) and \(\widetilde R\) are produced above. The previous
	estimates show that for suitable \(M\) large enough (depending on the data)
	and for \(T\) sufficiently small the image \(\mathcal G(\mathcal B_M)\subset\mathcal B_M\).
	Indeed the right-hand sides of \eqref{est:C_from_YR}, \eqref{est:Y_from_C},
	and \eqref{est:R_from_v} are bounded by constants depending on \(M\) and
	the data. If we choose \(M\) larger than those constants and then pick \(T\) small
	enough so exponentials and $T$ factors still do not break the bounds.\\
	It remains to  get the contraction property. Indeed, let \((\YY^1,R^1),(\YY^2,R^2)\in\mathcal B_M\) and denote the corresponding
	objects by \(\CC^1,\CC^2\), \(v^1,v^2\), \(\widetilde \YY^1,\widetilde \YY^2\),
	\(\widetilde R^1,\widetilde R^2\). We estimate the difference
	\(\delta := \|(\YY^1,R^1)-(\YY^2,R^2)\|_{\mathcal X_T}\).\\
	From the Lemma \ref{lem:main-final} and since \(\boldsymbol{H}\) is  Lipschitz 
	\begin{equation}\label{est:C_lip}
		\|\CC^1-\CC^2\|_{W^{2,1}_p}\le L_C(M)\big(\|\YY^1-\YY^2\|_{C(\OMV\times[0,T])} + \|R^1-R^2\|_{C([0,T])}\big),
	\end{equation}
	hence by embedding,
	\[
	\|\CC^1-\CC^2\|_{C(\OMV\times[0,T])}\le L_{C,emb}(M)\,\delta.
	\]
	Then using the definition of \(v\) and since and $g$ is also Lipschitz we have 
	\[
	\|v^1-v^2\|_{C(\OMV\times[0,T])}
	\le T\cdot L_v(M)\,\delta,
	\]
	where the factor \(T\) appears after integrating the difference of \(g(\YY,\CC)\)
	in time.
	Thus for small \(T\) the mapping to the velocity is contractive in the sense
	that its Lipschitz constant contains a factor \(T\).\\
	Next, compare the Duhamel representations of \(\widetilde \YY^1\) and
	\(\widetilde \YY^2\) and use Lipschitzness of \(\boldsymbol{F}\) and the estimate for
	\(\|v^1-v^2\|\). One finds
	\begin{equation}\label{est:Y_lip}
		\|\widetilde \YY^1-\widetilde \YY^2\|_{C^1}
		\le T\cdot L_Y(M)\,\delta,
	\end{equation}
	again with a small factor \(T\) (the transport Duhamel integral gives
	explicitly a factor \(t\) which is controlled by \(T\)).\\
	Using the Gronwall lemma and the standard ODE theory we get that
	\begin{equation}\label{est:R_lip}
		\|\widetilde R^1-\widetilde R^2\|_{C^1}
		\le T\cdot L_R(M)\,\delta.
	\end{equation}
	Combining \eqref{est:Y_lip}--\eqref{est:R_lip} we obtain, for a constant
	$L_{\mathcal G}(M)$,
	\[
	\|\mathcal G(\YY^1,R^1)-\mathcal G(\YY^2,R^2)\|_{\mathcal X_T}
	\le L_{\mathcal G}(M)\,T\,\delta.
	\]
	Hence choosing \(T\) so small that \(L_{\mathcal G}(M)T<1\) we see that
	\(\mathcal G\) is a contraction on \(\mathcal B_M\).\\
	By Banach fixed point theorem there exists a unique fixed point
	\((\YY,R)\in\mathcal B_M\) with \(\mathcal G(\YY,R)=(\YY,R)\). The parabolic solve
	with this \((\YY,R)\) gives the corresponding \(\CC\). Thus the triple
	\((\CC,\YY,R)\) is a solution of the coupled system on \([0,T]\). Uniqueness
	in the class \(\CC\in W^{2,1}_p,\ \YY\in C^1,\ R\in C^1\) follows from the
	contraction argument.\\
	Combining  the estimates from the parabolic Lemma \eqref{est:C_from_YR},
	the transport Lemma \eqref{est:Y_from_C}, and the ODE estimate
	\eqref{est:R_from_v}yield bounds of the form
	\[
	\|\CC\|_{W^{2,1}_p}\le A_1 + A_2\|\YY\|_{C^1} + A_3\|R\|_{C^1},
	\qquad
	\|\YY\|_{C^1}\le B_1 + B_2 t\sup_{[0,t]}\|\boldsymbol{F}\|_{C},
	\]
	\[
	\|R\|_{C^1}\le N_1 + N_2\|v(1,\cdot)\|_{C([0,t])},
	\]
	for constants depending on the data and Lipschitz constants. Iterating these
	inequalities and using Gronwall we obtain the stated exponential-type bound
	\eqref{eq:main_estimate} (with constants \(K_1,K_2\) depending only on the
	data and the local ball radius \(M\)). This is the required \textit{a priori} estimate on the local interval.\\
	Finally, standard continuation arguments apply: the solution can be extended
	past \(T^\ast\) as long as the norms appearing in
	\(\|\CC\|_{W^{2,1}_p}+\|\YY\|_{C^1}+\|R\|_{C^1}\) remain finite. This finishes
	the proof.
\end{pre}

\section{Global existence via invariance regions and energy estimates}
This section is dedicated to extend the results iobtained in the previous sections for any  $t\ge 0.$ To achieve this we need the following assumptions which extend the previous ones. Let us define $a(z,t)=zv(z,t).$ We have the following theorem
\begin{theo}\label{thm:global}
	Let $T_1>0$ and consider the coupled system \eqref{EqBiomassVersion2} --\eqref{EqThicknessVersion2} on $(z,t)\in\OMV\times[0,T_1].$\\
	Assume the following structural hypotheses hold.
	\begin{enumerate}
		\item $a\in C(\OMV\times[0,T_1]),$  $a(1,t)=a(0,t)=0$ and there exists $M_a>0$ such that
		$\|a\|_{L^\infty}\le M_a$. 
		
		\item (Diffusion) $D_j>0$ for all $j$.
		
		\item (Lipschitzness) $f=(f_i)_{i=1}^n:\mathbb R^n\times\mathbb R^m\to\mathbb R^n$
		and $h=(h_j)_{j=1}^m:\mathbb R^n\times\mathbb R^m\to\mathbb R^m$ are locally
		Lipschitz and satisfy the quasi-positivity property:
		\[
		f_i(\YY,\CC)\ge 0\ \text{ whenever }\YY\ge \boldsymbol{0},\CC\ge \boldsymbol{0},\]
		\[ h_j(\YY,\CC)\ge \boldsymbol{0}\ \text{ whenever }\YY\ge \boldsymbol{0},\CC\ge \boldsymbol{0}.
		\]
		\item (Dissipative energy inequality) There exist constants $\alpha>0$,
		$\beta\ge 0$ and $M_0\ge0$ such that for all $Y\in\mathbb R^n$, $C\in\mathbb R^m$
		and all $t\ge 0$,
		\begin{equation}\label{eq:dissipative}
			\sum_{i=1}^n \mu_i f_i(\YY,\CC) Y_i + \sum_{j=1}^m \nu_j h_j(\YY,\CC) C_j
			\le -\alpha\big(\|\YY\|_2^2 + \|\CC\|_2^2\big) + \beta + M_0,
		\end{equation}
		for some fixed positive weights $\mu_i,\nu_j>0$ (one may take $\mu_i=\nu_j=1$).
	\end{enumerate}
	Then the unique local classical solution given by the local existence theorem \ref{thm:main_existence1} extends to a \emph{global} classical
	solution on $\OMV\times[0,\infty)$ with
	\[
	\YY\in C^1(\OMV\times[0,\infty))^n,~
	\CC\in W^{2,1}_{p,\mathrm{loc}}(\OMV\times[0,\infty))^m,~v\in C^1(\OMV\times[0,\infty))^n,
	R\in C^1([0,\infty)).
	\]
	Moreover, the solution remains nonnegative and satisfies the uniform in time
	energy bound
	\begin{equation}\label{eq:global_bound}
		\|\YY(\cdot,t)\|_{L^2}^2 + \|\CC(\cdot,t)\|_{L^2}^2 + \|v(\cdot,t)\|_{L^2}^2+ R^2(t)\le C_\infty
		\qquad\text{for all } t\ge 0,
	\end{equation}
	for some constant $C_\infty$ depending only on the initial data and model
	constants.
\end{theo}
This inequality \eqref{eq:global_bound} is a dissipation condition which ensures  that reaction terms
do not allow arbitrary growth and provide a restoring effect for large norms. To prove this system we  use the local existence result and the continuation
criterion then  derive \textit{a priori} energy estimates which are uniform on finite
time intervals and in fact globally bounded. We conclude the proof  by using  these estimates to
prevent the norms that control continuation from blowing up; hence extend the
solution globally.
\begin{pre} By the local well-posedness result in Theorem \ref{thm:main_existence1}) we only  need to obtain \textit{a priori} bounds preventing blow-up of the estimates obtained earlier.
	In addition, since each $f_i$ and $h_j$ are quasi-positive (hypothesis 3)) and the
	transport and diffusion operators preserve sign under the given boundary
	conditions, standard comparison-principle arguments yield
	\[
	Y_i(z,t)\ge 0,\qquad C_j(z,t)\ge 0\qquad\text{for all }(z,t)\in[0,1]\times[0,T_{\mathrm{loc}}].
	\]
	Indeed, if at some time a component would become negative, the quasi-positivity
	of the source forbids a decrease through the source term; transport and
	diffusion alone preserve nonnegativity when initial or boundary data are
	nonnegative. From the local solution this property holds on $[0,T_{\mathrm{loc}}]$.
	Thus the region \(\{\YY\ge0,\CC\ge0\}\) is invariant.\\
	To get the energy estimates we multiply each biomass equation \eqref{EqBiomassVersion2} by $\mu_i Y_i$ and integrate
	over $z\in(0,1)$; multiply each substrate equation \eqref{ParabolicVersion2} by
	$\nu_j C_j$ and integrate, then sum the results (weights $\mu_i,\nu_j>0$ are
	those appearing in \eqref{eq:dissipative}). Using integration by parts for
	the diffusion terms 
	we obtain
	\begin{align}\label{eq:energy_derivative}
		\frac{d}{dt}\Big\{\tfrac12\|\YY(\cdot,t)\|_{L^2_\mu}^2 + \tfrac12\|\CC(\cdot,t)\|_{L^2_\nu}^2\Big\}
		+\; \sum_{j=1}^m \nu_j D_j \int_0^1 |\partial_z C_j|^2\,dz \nonumber\\
		=\; R^2(t)\int_0^1\Big(\sum_{i=1}^n \mu_i f_i(\YY,\CC)Y_i + \sum_{j=1}^m \nu_j h_j(\YY,\CC) C_j\Big)\,dz \nonumber\\
		\qquad + \int_0^1\Big(zv(1,t)\partial_z C_j C_j+ zv(1,t)\partial_z Y_i Y_i\Big)dz,
	\end{align}
	where $\|Y\|_{L^2_\mu}^2=\sum_i\mu_i\int_0^1 Y_i^2\,dz$ and similarly for $C$.
	Using the hhypothesis 1, the last two integral vanish. For
	generality we assume boundary contributions are controlled (they are either
	zero under common boundary conditions or can be estimated by the $L^2$-norms of $\YY,\CC$). \\
	Next, we apply \eqref{eq:dissipative} pointwise in $(z,t)$ and integrate over $z$
	to estimate the remainder of right-hand side of \eqref{eq:energy_derivative}:
	\[
	R^2(t)\int_0^1\Big(\sum_{i}\mu_i f_i Y_i + \sum_{j}\nu_j h_j C_j\Big)\,dz
	\le R^2(t)\int_0^1\big(-\alpha(\|Y\|_2^2+\|C\|_2^2)+\beta+M_0\big)\,dz.
	\]
	Because $\int_0^1(\|\YY\|_2^2+\|\CC\|_2^2)\,dz=\|\YY\|_{L^2}^2+\|\CC\|_{L^2}^2$, we get
	(omitting the boundary term for clarity)
	\begin{equation}\label{eq:energy_ineq}
		\frac{d}{dt}\mathcal{E}(t) + \kappa\|\partial_z \CC(\cdot,t)\|_{L^2}^2
		\le -\alpha R^2(t)\big(\|\YY(\cdot,t)\|_{L^2}^2+\|\CC(\cdot,t)\|_{L^2}^2\big) + R^2(t)(\beta+M_0),
	\end{equation}
	where we have set the energy
	\[
	\mathcal{E}(t):=\tfrac12\big(\|\YY(\cdot,t)\|_{L^2_\mu}^2 + \|\CC(\cdot,t)\|_{L^2_\nu}^2\big)
	\]
	and $\kappa=\min_j \nu_j D_j>0$.
	To bound  R we recall that it satisfies the ODE \eqref{EqThicknessVersion2}. One notices that the ODE
	contains a strong damping term $-\lambda R^4$ which prevents runaway of
	$R$ for large values. More precisely if we consider the function $Q(t):=R^2(t)$, then we have 
	\[
	\dot Q(t)=2R(t)\dot R(t)=2R^3(t)v(t)-2\lambda R^5(t)=2Q^{3/2}v(t)-2\lambda Q^{5/2}.
	\]
	From this we get a differential inequality that shows $R(t)$ cannot blow up
	in finite time because for large $Q$ the $-2\lambda Q^{5/2}$ term dominates and forces
	a decrease. More directly, using standard comparison with the autonomous ODE
	$\dot R=-\tfrac12\lambda R^4$ for sufficiently large $R$ and since $v$ is bounded, we further infer that there exists $R_{\max}$ depending on $R_0,\lambda$
	and $\|v\|_\infty$ such that
	\[
	R(t)\le R_{\max}\qquad\text{for all }t\ge0.
	\]
	In particular $R^2(t)$ is globally bounded. Notably there exists $M_R>0$ such that
	$R^2(t)\le M_R$ for all $t\ge0$ makes sense.\\
	Using the bound $R^2(t)\le M_R$ in \eqref{eq:energy_ineq} yields
	\[
	\frac{d}{dt}\mathcal{E}(t)
	\le -\alpha M_R \big(\|\YY\|_{L^2}^2+\|\CC\|_{L^2}^2\big) + M_R(\beta+M_0).
	\]
	Dropping the negative first term on the right gives the crude linear bound
	\[
	\frac{d}{dt}\mathcal{E}(t)\le M_R(\beta+M_0).
	\]
	Integrating from $0$ to $t$,
	\[
	\mathcal{E}(t)\le \mathcal{E}(0) + M_R(\beta+M_0)\,t.
	\]
	This alone gives at most linear growth. To obtain a uniform-in-time bound we
	use the full dissipative structure and rewrite \eqref{eq:energy_ineq} as
	\[
	\frac{d}{dt}\mathcal{E}(t) \le -\alpha M_R\big(2\mathcal{E}(t)/C^*\big) + M_R(\beta+M_0),
	\]
	where $C^*>0$ is a constant comparing the weighted and usual $L^2$-norms:
	$\|\YY\|_{L^2_\mu}^2 + \|\CC\|_{L^2_\nu}^2 \ge C^*(\|\YY\|_{L^2}^2 + \|\CC\|_{L^2}^2)$.
	Thus
	\[
	\frac{d}{dt}\mathcal{E}(t) \le -\gamma \mathcal{E}(t) + M_R(\beta+M_0),
	\qquad\text{with }\gamma=\frac{2\alpha M_R}{C^*}>0.
	\]
	By the standard linear Gronwall  inequality for this linear
	dissipative ODE we obtain
	\[
	\mathcal{E}(t) \le e^{-\gamma t}\mathcal{E}(0) + \frac{M_R(\beta+M_0)}{\gamma}
	\le \mathcal{E}(0) + \frac{M_R(\beta+M_0)}{\gamma}=:E_\infty.
	\]
	Hence $\mathcal{E}(t)$ is uniformly bounded for all $t\ge 0$.
	
	Consequently there exists $C_\infty>0$ so that
	\[
	\|\YY(\cdot,t)\|_{L^2}^2 + \|\CC(\cdot,t)\|_{L^2}^2 \le C_\infty,\qquad\forall t\ge0,
	\]
	which, together with the bound for $R(t)$, yields \eqref{eq:global_bound}. Notice the velocity bound can also be deduced by using the fact that $g$ is sum of terms from $f_i$ (which are bounded) and boundeness of $R.$\\
	We now upgrade the $L^2$-bounds to the norms used in the continuation
	criterion. For the parabolic variables $\CC$: since the right-hand side
	$R^2(t)h_j(\YY,\CC)$ is uniformly bounded in $L^2(\OM)$ by the $L^2$-bounds on
	$\YY,\CC$ and the uniform bound on $R$, classical $L^2$-parabolic regularity
	estimates yield that for any
	finite time $T>0$,
	\[
	\|\CC\|_{W^{2,1}_p(\OM\times(0,T))}\le C_G(T),
	\]
	for some finite constant $C_G(T)$ depending only on the initial data and the
	bounds obtained above (and on $T$). In particular for each finite $T$ the
	$\CC$-norms used in the continuation criterion remain finite.\\
	For the transport variables $\YY$: the transport equation admits uniform
	$L^\infty$–bounds because the right-hand side $R^2 f_i(\YY,\CC)$ is controlled
	by the $L^2$-bounds on $Y,C$ and the Lipschitz structure of $f_i$ together
	with invariant region arguments. More precisely, using the Duhamel integral
	formula along characteristics and Gronwall-type arguments, one deduces for
	each finite $T$ a bound
	\[
	\|\YY\|_{C^1(\OMV\times[0,T])} \le C_Y(T),
	\]
	again depending only on the \textit{a priori} quantities already bounded. 
	Finally, $R(t)$ is bounded and continuous for all time as shown earlier. Therefore none of the norms that control the local continuation criterion blow up in finite time.\\
	Since all continuation norms remain bounded for all finite times, the local
	solution can be extended step by step to a global solution on $[0,\infty)$.
	The uniform-in-time \textit{a priori} bounds \eqref{eq:global_bound} prevent finite-time
	blow-up and guarantee classical regularity persists for all times.
\end{pre}

\newpage
\section{Statements and Declarations}

\textbf{Competing Interests:} The author declares no competing interests.

\end{document}